\documentclass[psamsfonts]{amsart}

\newtheorem{theorem}{Theorem}[section]
\newtheorem{lemma}[theorem]{Lemma}
\newtheorem{proposition}[theorem]{Proposition}
\newtheorem{corollary}[theorem]{Corollary}

\theoremstyle{definition}
\newtheorem{definition}[theorem]{Definition}
\newtheorem{example}[theorem]{Example}

\newtheorem{remark}[theorem]{Remark}

\newtheorem{situation}[theorem]{Situation}

\theoremstyle{remark}

\numberwithin{equation}{section}

\newcommand{\NN}{\mathbb{N}}
\newcommand{\ZZ}{\mathbb{Z}}

\newcommand{\RR}{\mathbb{R}}

\newcommand{\PP}{\mathbb{P}}

\newcommand  {\shE}     {\mathcal{E}}
\newcommand  {\shF}     {\mathcal{F}}
\newcommand  {\shG}     {\mathcal{G}}
\newcommand  {\shH}     {\mathcal{H}}

\newcommand  {\shM}     {\mathcal{M}}

\newcommand  {\shN}     {\mathcal{N}}
\newcommand  {\shL}     {\mathcal{L}}
\newcommand  {\shR}     {\mathcal{R}}
\newcommand  {\shS}     {\mathcal{S}}
\newcommand  {\shT}     {\mathcal{T}}

\newcommand  {\shQ}     {\mathcal{Q}}

\newcommand  {\foa}     {\mathfrak{a}}

\newcommand  {\fop}     {\mathfrak{p}}
\newcommand  {\foq}     {\mathfrak{q}}

\newcommand  {\Char}    {\operatorname{char}}

\newcommand  {\Det}    {\operatorname{Det}}

\newcommand  {\dual}    {\vee}

\newcommand  {\Ext}     {\operatorname{Ext}}

\newcommand  {\Hom}     {\operatorname{Hom}}

\newcommand  {\lra}     {\longrightarrow}

\newcommand  {\mult}    {\operatorname{mult}}

\renewcommand{\O}       {\mathcal{O}}

\newcommand  {\Proj}    {\operatorname{Proj}}

\newcommand  {\ra}      {\rightarrow}

\newcommand  {\rk}    {\operatorname{rk}}

\newcommand  {\Spec}    {\operatorname{Spec}}

\newcommand {\mindeg} {\rho}

\def\mydate{\number\day\space\ifcase\month \or January\or February\or March\or April\or May\or
June\or July\or August\or September\or October\or November\or
December\fi \space\number\year}

\usepackage{amscd}
\usepackage{amssymb}

\begin{document}

\title[Slopes and tight closure]
{Slopes of vector bundles on projective curves and applications
to tight closure problems}


\author[Holger Brenner]{Holger Brenner}
\address{Mathematische Fakult\"at, Ruhr-Universit\"at Bochum, 
               44780 Bochum, Germany}
\email{brenner@cobra.ruhr-uni-bochum.de}





\begin{abstract}
We study different notions of slope of a vector bundle
over a smooth projective curve with respect to ampleness and affineness
in order to apply this to tight closure problems.
This method gives new degree estimates from above and from below
for the tight closure of a homogeneous $R_+$-primary ideal
in a two-dimensional normal standard-graded algebra $R$ in terms of the minimal
and the maximal slope of the sheaf of relations
for some ideal generators. If moreover this sheaf of relations is semistable,
then both degree estimates coincide and we get a vanishing type theorem.
\end{abstract}

\maketitle

\section*{Introduction}

In this paper we continue the study of tight closure problems
started in \cite{brennertightproj}
of a two-dimensional normal standard-graded $K$-algebra $R$ in terms of the corresponding
projective bundles and subbundles on the corresponding
projective curve $Y= \Proj\, R$.
A system of $R_+$-primary homogeneous elements $f_1, \ldots ,f_n \in R$
leads to the locally free sheaf of relations $\shR$ on $Y$, and a further
homogeneous element yields an extension
$\shR \ra \shR'$, which itself gives the projective subbundle
$\PP(\shF) \subset \PP(\shF')$, where $\shF= \shR^\dual$
and $\shF' = {\shR'}^\dual$.

The fundamental observation of \cite{brennertightproj}
is that $f_0$ belongs to the tight closure of
the ideal $(f_1, \ldots ,f_n)$, that is $f_0 \in (f_1, \ldots ,f_n)^*$ holds
if and only if the open subset $\PP(\shF') -\PP(\shF)$ is not an affine scheme.
This link rests upon the reinterpretation of tight closure as
solid closure in positive characteristic, see \cite{hochstersolid}
for this notion
and \cite{hunekeapplication} and \cite{hunekeparameter}
for background on the theory of tight closure
(in characteristic zero we work throughout this paper with the
notion of solid closure).
This gives a powerful geometric tool to study tight closure problems,
see also \cite{brennertightproj} for the general setting
and \cite{brennertightelliptic} for the proof that tight closure
and plus closure are the same for $R_+$-primary ideals
in a normal coordinate ring of an elliptic curve over a
field of positive characteristic.

The theme of this paper is the minimal and the maximal slope
of a locally free sheaf $\shG$ on a smooth projective curve $Y$
and how these invariants and the corresponding properties
like semistability and ampleness are related to the properties
which are of interest from the tight closure point of view.
We recall and extend the necessary definitions and facts in
section \ref{sectionslope}. In section \ref{sectionample}
we consider criteria for ample bundles.
Our main criterion, which is essentially due to Barton
(\cite[Theorem 2.1]{barton}), is that $\shG$ is ample if and only if
$\bar{\mu}_{\min} (\shG) >0$ (Theorem \ref{amplekrit}), where
$\bar{\mu}_{\min} (\shG)$ is a variant of the minimal slope
which takes also into account the behavior under finite mappings $Z \ra Y$.

In section \ref{sectionaffine} we study sufficient conditions
for the complement $\PP(\shG') -\PP(\shG)$ to be affine,
where $0 \ra \O_Y \ra \shG' \ra \shG \ra 0$
is an extension given by a cohomology class $c \in H^1(Y, \shG^\dual)$.
From the ampleness criterion it follows in characteristic $0$ easily that
for $\mu (\shG) >0$ and $c \neq 0$ the complement
$\PP(\shG') -\PP(\shG)$ is affine (Theorem \ref{ampleaffin}).
This is also true if there exists a sheaf homomorphism
$\varphi: \shG^\dual \ra \shT$ such that
$0 \neq \varphi(c) \in H^1(Y, \shT)$ and $\shT$ is semistable of negative
slope (Theorem \ref{slopekritaffin}).
Our main result of section \ref{sectionnonaffine} is that
$\PP(\shG') -\PP(\shG)$ is not affine under the condition
that $\bar{\mu}_{\max}(\shG) \leq 0$ (Theorem \ref{slopemaxkrit}).
The sections \ref{sectionslope} - \ref{sectionnonaffine} may be read without
any knowledge of tight closure theory.

In section \ref{sectiongraded}
we recall briefly how vector bundles arise from tight closure problems
referring to \cite{brennertightproj} for details and proofs.
In sections \ref{sectioninclusion} - \ref{examples}
we apply our results to tight closure problems.
The slope invariants $\bar{\mu}_{\max}$ and $\bar{\mu}_{\min}$
of the relation bundle corresponding
to ideal generators $f_1, \ldots ,f_n$ of a $R_+$-primary ideal
give important new degree estimates for $(f_1, \ldots ,f_n)^*$.

We show that the condition
$\deg \, (f_0) \geq \bar{\mu}_{\max} (f_1, \ldots ,f_n)/\deg \, (\O_Y(1))$
forces $f_0 \in (f_1, \ldots ,f_n)^*$ (Theorem \ref{maxin}).
From this we derive that
$R_m \subseteq (f_1, \ldots ,f_n)^*$, whenever
$m$ is greater or equal the sum of the two biggest degrees of the $f_i$
(Corollary \ref{inclusionbound}),
which improves slightly the bound $2 \max d_i$ of Smith given
in \cite[Proposition 3.1]{smithgraded}.

In the other direction we prove
(if the characteristic of the field is $0$ or $p  \gg  0$) that if
$\deg \, (f_0) < \bar{\mu}_{\rm min} (f_1, \ldots ,f_n)/\deg \, (\O_Y(1))$,
then $f_0 \in (f_1, \ldots, f_n)^*$
is only possible if already 
$f_0 \in (f_1, \ldots ,f_n)$ holds (Theorem \ref{minex}).
If the sheaf of relations splits into invertible sheaves, then the
two bounds are easy to compute as the minimum (maximum) of the degrees
of the summands.
Moreover, in this splitting situation we can give a numerical criterion
for tight closure (Theorem \ref{splitting}).

In section \ref{sectionvanishing} we study the situation where
the minimal and the maximal slope coincide.
In this case the sheaf of relations is semistable, and we get the equality
$(f_1, \ldots ,f_n)^*=(f_1, \ldots ,f_n) + R_{\geq k}$,
where $k= \lceil \frac{d_1 + \ldots +d_n}{n-1} \rceil $
(Theorems \ref{semistablevanishing} and \ref{semistablevanishingp}). 
This is a (so-called) vanishing type theorem 
in the sense of \cite{hunekesmithkodaira}
and generalizes the vanishing theorem in the parameter case ($n=2$),
where the sheaf of relations is invertible, hence of course semistable.

In the last section \ref{examples} we concentrate on the case
$n=3$. We give bounds for $\mu_{\rm max}$ and $\mu_{\rm min}$
deriving from known results about the $e$-invariant of a ruled surface
(Theorem \ref{genusbound}).
Under the condition that the sheaf of relations is indecomposable
we obtain degree bounds for inclusion and exclusion
which are quite near to $k=(d_1+d_2+d_3)/2$, the difference is at most
$(g-1)/\deg (\O_Y(1))$, where $g$ is the genus and of the curve
(Corollary \ref{genusboundtight}).


\section{Slope of bundles}
\label{sectionslope}

In this section we recall the notion of the slope of bundles
and various related concepts like semistable sheaves
and the Harder-Narasimhan filtration which we will need in the sequel,
our main references are \cite{hardernarasimhan}, \cite{huybrechtslehn},
\cite[Ch. 6.4]{lazarsfeld}, \cite{miyaokachern} and \cite{seshadrifibre}.

Let $\shE$ denote a locally free sheaf on a smooth projective
curve $Y$ over an algebraically closed field $K$.
The degree of $\shE$ is defined by $\deg \, (\shE) = \deg \, (\bigwedge^r (\shE))$,
where $r$ is the rank of $\shE$.
If $\PP(\shE)= \Proj \, \oplus_k S^k(\shE)$
is the corresponding projective bundle of
dimension $r$ and if $\xi$ denotes the divisor class corresponding
to the relatively very ample invertible sheaf
$\O_{\PP(\shE)}(1)$, then also $\deg \, (\shE)= \xi^r$
equals the top self intersection number, see \cite[Lemma 6.4.10]{lazarsfeld}.
The number $\deg \, (\shE)/ r!$ is also
the coefficient of $k^r$ in the Euler-Hilbert polynomial
$\chi( \O_{\PP(\shE)}(k))$, which equals $\chi( S^k(\shE))$.

The {\em slope} of a locally free sheaf $\shE$ is defined by
$\mu (\shE)= \deg \, (\shE) / \rk (\shE)$.
A locally free sheaf $\shE$ is called {\em semistable}, if for every
locally free quotient sheaf $\shE \ra \shQ \ra 0$ the inequality
$\mu (\shQ) \geq \mu(\shE)$ holds.
This is equivalent to the property that for every locally free
subsheaf $\shT \subseteq \shE$ the inequality
$\mu (\shT) \leq \mu (\shE)$ holds.

Every locally free sheaf $\shE$
has a unique {\em Harder-Narasimhan Filtration}.
This is a filtration of locally free subsheaves
$$ 0=\shE_0 \subset \shE_1 \subset \ldots \subset \shE_s = \shE $$
such that $\shE_i/\shE_{i-1}$ is semistable for every $i=1, \ldots, s$.
$\shE_1$ is called the {\em maximal destabilizing subsheaf}.
The slopes of these semistable quotients
form a decreasing chain
$\mu_1 > \ldots > \mu_s$.
$\mu_{\rm min}(\shE)= \mu_s= \mu(\shE/\shE_{s-1})$ is called the
{\em minimal slope}
and $\mu_{\max}(\shE) =\mu_1 (\shE)$ is called the {\em maximal slope}.
This is the same as
$\mu_{\rm min}(\shE)=\min\{\mu(\shQ):\, \shE \ra \shQ \ra 0 \}$.
For the dual sheaf we have
$\mu_{\rm max} (\shE^\dual )= -\mu_{\rm min}(\shE)$.
If $\mu_{\min} (\shE) > \mu_{\max}(\shF)$, then $\Hom(\shE,\shF)=0$.
In particular, if $\mu_{\rm max} (\shF)<0$,
then $\Gamma(Y, \shF)=0$.

If $\varphi: Z \ra Y$ is a finite $K$-morphism between smooth projective
curves over the algebraically closed field $K$, then
$\mu(\varphi^*(\shE))= \deg \, (\varphi ) \mu (\shE)$.
If $\varphi $ is separable, then
the pull-back $ \varphi^*(\shE)$ of a semistable sheaf
$\shE$ on $Y$ is again semistable, see \cite[Proposition 3.2]{miyaokachern}.
Hence in the separable case
the Harder-Narasimhan filtration of $\varphi^*(\shE)$ is just the
pull-back of the Harder-Narasimhan filtration of $\shE$,
and also $\mu_{\max}$ and $\mu_{\min}$ transform in the same way
as $\mu$ does.

In the non-separable case this is not true at all
and the notion of semistability needs to be refined.
A locally free sheaf $\shE$
on $Y$ is called {\em strongly semistable},
if for every finite $K$-morphism $\varphi:Z \ra Y$ the pull-back
$\varphi^*(\shE)$ is again semistable. In characteristic zero,
this is the same as being semistable,
and in positive characteristic it is given by the property that
the pull-back under every $K$-linear Frobenius morphism is semistable,
see \cite[Proposition 5.1]{miyaokachern}.
This difficulty in positive characteristic
is one motivation for the following definition.

\begin{definition}
Let $Y$ denote a smooth projective curve over an algebraically closed
field and let
$\shE$ denote a locally free sheaf. Then we define
$$
\bar{\mu}_{\max} (\shE)=
{\rm sup} \{  \frac{\mu_{\rm max} (\varphi^* \shE)}{\deg \, (\varphi) }| \,\,
\varphi: Z \ra Y \mbox{ finite dominant $K$-morphism }      \}
$$
and
$$
\bar{\mu}_{\rm min} (\shE)=
{\rm inf} \{  \frac{\mu_{\rm min} (\varphi^* \shE)}{\deg \, (\varphi)} |\, \,
\varphi: Z \ra Y \mbox{ finite dominant $K$-morphism }      \} \, .
$$
\end{definition}

\begin{remark}
It is enough to consider in the previous definition
only $K$-linear Frobenius morphisms, since every morphism factors
through a separable map and a Frobenius and the maximal and minimal
slope behave well with respect to separable morphisms.

We will see in remark \ref{ampleremark} that $\bar{\mu}_{\min} (\shE)$
is bounded from below,
hence these numbers exist, but it is not clear whether they are obtained.
An equivalent question is whether one may find
a sufficiently high Frobenius power
such that the Harder-Narasimhan filtration of $\varphi^*(\shE)$ consists
of strongly semistable quotients.
\end{remark}

We will also need the following definition,
compare with \cite{ionescutoma}, \cite{langeregel}
and \cite{mukaisakai}
for related invariants of $\shE$ and $\PP(\shE)$.

\begin{definition}
Let $\shE$ denote
a locally free sheaf on a smooth projective curve,
let $1 \leq s \leq \rk \, (\shE)$. We set
$$\mindeg_s \, (\shE) = \min \{ \deg \, (\shQ): \, \shE \ra \shQ \ra 0, \, \,\,
\shQ \mbox{ is locally free and } \rk \, (\shQ) =s \} \, .$$
\end{definition}

\section{Ampleness criteria for vector bundles over projective curves}
\label{sectionample}

Recall that a locally free sheaf $\shG$ on a scheme $Y$ is called {\em ample}
if the invertible sheaf $\O_{\PP(\shG)}(1)$ on the projective bundle
$\PP(\shG) = \Proj \, \oplus_n S^n(\shG)$ is ample.
In characteristic zero we have the following
linear criterion of Hartshorne-Miyaoka for ample sheaves over a curve.

\begin{theorem}
\label{amplekritnull}
Let $Y$ denote a smooth projective curve over an algebraically closed
field of characteristic $0$.
Let $\shG$ denote a locally free sheaf on $Y$ and
$\PP(\shG)$ the corresponding projective bundle,
and let $\xi$ denote the divisor class corresponding to $\O_{\PP(\shG)}(1)$.
Then the following are equivalent.

\renewcommand{\labelenumi}{(\roman{enumi})}
\begin{enumerate}

\item
The sheaf $\shG$ is ample.

\item
For every projective subbundle
$\PP(\shQ) \subseteq \PP(\shG)$
of dimension $s$ we have
$\PP(\shQ) . \xi^s >0 $.

\item
For every locally free quotient sheaf $\shG \ra \shQ \ra 0$, $\shQ \neq 0$,
we have $\deg \, (\shQ) > 0$.

\item
The minimal slope is $\mu_{\rm min}(\shG) >0$.
\end{enumerate}

\end{theorem}
\proof
The equivalence of (ii),(iii) and (iv) is clear
since $\PP(\shQ). \xi^s =(\xi| _{\PP(\shQ)})^s = \deg \, (\shQ)$.
(i) $\Rightarrow$ (ii) follows from the Nakai-criterion,
for the other direction see
\cite[Theorem 2.4]{hartshorneamplecurve},
\cite[Corollary 3.5]{miyaokachern} or \cite[Theorem 6.4.15]{lazarsfeld}.
\qed

\bigskip
The minimal degree $\mindeg_1(\shG)$ of a quotient invertible sheaf of
$\shG$ must fulfill a stronger condition to guarantee that
$\shG$ is ample.

\begin{corollary}
Let $Y$ denote a smooth projective curve of genus $g$
over an algebraically closed field of characteristic $0$.
Let $\shG$ denote a locally free sheaf of rank $r$ on $Y$ and
suppose that
$\mindeg_1(\shG) > \frac{r-1}{r} g$.
Then $\shG$ is ample.
\end{corollary}
\proof
This is deduced in \cite[Proposition 2]{ionescutoma}
from the theorem.
\qed

\bigskip
The conditions (ii)-(iv) in Theorem \ref{amplekritnull}
are in positive characteristic not sufficient
for ampleness.
This is due to the fact that the pull back of a semistable sheaf
under a non separable morphism need not be semistable anymore and the
minimal slope may drop.
It is also related to the failure of the vanishing theorem
in tight closure theory
in small positive characteristic. The following criterion is valid
in every characteristic and is essentially due to Barton
(see \cite[Theorem 2.1]{barton}).

\begin{theorem}
\label{amplekrit}
Let $Y$ denote a smooth projective curve over an algebraically closed
field $K$ of characteristic $p \geq 0$.
Let $\shG$ denote a locally free sheaf of rank $r$ on $Y$.
Then the following are equivalent.

\renewcommand{\labelenumi}{(\roman{enumi})}
\begin{enumerate}

\item
The sheaf $\shG$ is ample.

\item
$\bar{\mu}_{\rm min}(\shG) >0 $.

\item
There exists $\epsilon >0$ such that for every finite $K$-morphism
$\varphi :Y' \ra Y$ and every invertible quotient sheaf
$\varphi^*(\shG) \ra \shM \ra 0$
the inequality $\frac{\deg \, (\shM)}{ \deg \, (\varphi)} \geq \epsilon > 0$
holds.

\end{enumerate}

\end{theorem}

\proof
(i) $\Rightarrow $ (ii).
Suppose that
$\varphi: Y' \ra Y$ is finite, where $Y'$ is another
smooth projective curve,
and let $\varphi ^* (\shG) \ra \shH \ra 0$ be given,
$s= \rk (\shH)$.
We consider first the case $s=1$.
Then $\shM=\shH$ is an invertible sheaf on $Y'$ and
the surjection $\varphi ^* (\shG) \ra \shM \ra 0$
defines a section
$s: Y' \ra \PP(\varphi^*(\shG))$
due to the correspondence described in
\cite[Proposition 4.2.3]{EGAII}
and a curve $Z$ (its image)
in $\PP(\shG)$.
The map $Y' \ra Z$ factors through the normalization of $Z$,
hence $\shM$ is defined already on this normalization. Therefore
we may assume that $Y'$ is the normalization of $Z$.
The numerical class of the curve $Z$ in $\PP(\shG)$
can be written as
$a\xi^{r-1}+ b \xi^{r-2}.f$,
where $\xi$ is again the divisor class corresponding to $\O_{\PP(\shG)}(1)$,
$f$ is the class of a fiber $\PP(\shG) \ra Y$ and $a,b \in \ZZ$.
Furthermore $a$ equals the degree of $\varphi$.
Therefore we have
$$\deg \, (\shM) = Y' . s^*(\O_{\PP(\shG)} (1))
= ( a\xi^{r-1}+ b \xi^{r-2}.f) . \xi=a \xi^r +b
= a \deg \, (\shG) +b \, .$$
Hence
$\deg \, (\shM) / \deg \, (\varphi)  = \deg \, (\shG) + b/a $.
Since $\shG$ and hence $\xi$ is ample, there exists a number $\epsilon >0$
such that $ \xi.Z \geq \epsilon || Z||$ holds for every curve,
where $|| Z||$ denotes any norm on ${\rm Num} \otimes \RR$,
see \cite[Theorem 8.1]{haramp}. Hence
$$ \deg \, (\shG) +\frac{b}{a} = \frac{ \xi. Z}{a}
\geq \frac{\epsilon \sqrt{a^2+b^2}}{a } = \epsilon \sqrt{1+ (\frac{b}{a})^2}
\geq
\epsilon \, .
$$

Now consider the general case.
Since $\shG$ is ample, also its wedge product
$\bigwedge^s \shG$ is ample due to \cite[Corollary 2.6]{hartshorneample}.
The surjection
$\varphi^* (\shG) \ra \shH \ra 0$ yields a surjection
$\bigwedge^s(\varphi^* \shG) \ra \bigwedge^s \shH \ra 0$.
$\shM=\bigwedge^s \shH$ is invertible and since
$\bigwedge^s (\varphi^* \shG)$ is ample, there exists
an $\epsilon_s$ such that
$\deg \, (\bigwedge^s \shH)/ \deg \, (\varphi) \geq \epsilon_s >0$.
Then $\bar{\mu}_{\min} (\shG) \geq \min _s \epsilon_s/s >0$.

\smallskip
(ii) $\Rightarrow $ (iii) is a restriction to invertible
quotient sheaves.

\smallskip
(iii) $\Rightarrow $ (i).
Let $\xi$ denote the hypersection divisor corresponding to
$\O_{\PP(\shG)}(1)$.
Due to the ampleness criterion of Seshadri, see \cite[I\S7]{haramp},
it is enough to show that there exists an
$\epsilon >0$ such that
$\frac{ \xi.Z}{ m(Z)} \geq \epsilon >0$ holds for every
(effective) curve $Z$,
where $m(Z)= {\rm sup} \{ \mult_P(Z)\}$
is the maximal multiplicity of a point on $Z$.
So suppose that $Z$ is an irreducible curve in $\PP(\shG)$.
If $Z$ lies in a fiber $F \cong \PP^{r-1}$, then
$\xi . Z =\deg \, (Z) \geq m(Z)$.
Hence we may assume that $Z$ dominates the base.
Let $Z'$ be the normalization of $Z$, $i: Z' \ra  \PP(\shG)$ the corresponding
mapping and let $\varphi :Z' \ra Y$ be the composition.
Let $\varphi^* \shG \ra \shM \ra 0$ be the corresponding surjection
onto the invertible sheaf $\shM$.

The multiplicity $m(Z)$ is bounded above by $\deg \, (\varphi)$.
Therefore we have
$$ \frac{\xi.Z}{ m(Z)}
=\frac{\deg \, (\shM)}{m(Z)}
\geq \frac{\deg \, (\shM)}{ \deg \, (\varphi)}
\geq \epsilon > 0 \, .
$$
\qed

\begin{remark}
\label{ampleremark}
If $\shG$ is locally free on $Y$ and if $\O_Y(1)$ is an ample invertible sheaf
on $Y$,
then $\shG(n)= \shG \otimes \O_Y(n)$ corresponds to the invertible
sheaf $\O_{\PP(\shG)}(1) \otimes p^*(\O_Y(n))$ on $\PP(\shG)$,
see \cite[Proposition 4.1.4]{EGAII}.
Choosing $n$ high enough we may achieve that
$\shG(n)$ becomes ample. Since the slopes transform like
$\mu_{\rm min}( \shG \otimes \shL^k )
= \mu_{\rm min} (\shG) +k \deg \, (\shL) $ it follows that
$\bar{\mu}_{\min}(\shG)
= \bar{\mu}_{\rm min}(\shG \otimes \O_Y(k)) -k \deg \, (\O_Y(1))$
is bounded from below.
Dually it follows that $\bar{\mu}_{\max} (\shS)$ is bounded from above,
so both numbers exist (but it is not clear whether they are obtained).
\end{remark}

\begin{corollary}
\label{amplekritsemistable}
Let $Y$ denote a smooth projective curve over an algebraically closed
field of characteristic $p \geq 0$.
Let $\shG$ denote a locally free sheaf on $Y$.
Suppose that $\shG$ is strongly semistable.
Then $\shG$ is ample if and only if its degree is positive.
\end{corollary}
\proof
This follows directly from \ref{amplekrit}.
For another proof see \cite[Corollary 3.5 and \S 5]{miyaokachern}.
\qed

\medskip
The ampleness of a locally free sheaf has also the
following consequence on $\mu_{\max}$, which we will use
in section \ref{sectioninclusion}.

\begin{corollary}
\label{mumaxabschaetzung}
Let $\shG$ denote an ample locally free sheaf of rank
$r$ on a smooth projective curve $Y$.
Then we have the estimates (set $\mindeg_0 (\shG)=0$)
$$ \mu_{\rm max}(\shG) \leq
{\rm max}_{s=0, \ldots, r-1}
\frac{ \deg \, (\shG) - \mindeg_s (\shG)}{r-s} 
\leq \deg \, (\shG)
\, .$$
\end{corollary}

\proof
Let $\shT \subset \shG$ denote a locally free subsheaf
of positive rank with a short exact sequence
$ 0 \ra \shT \ra \shG \ra \shQ \ra 0$,
$s= \rk (\shQ)$, $s=0, \ldots, r-1 $.
Then
$$\mu (\shT)= \frac{\deg \, (\shT)}{r-s}
=\frac{ \deg \, (\shG) -\deg \, (\shQ)}{r-s}
\leq \frac{ \deg \, (\shG) -\mindeg_s\, (\shG)}{r-s}
$$
due to the definition of
$\mindeg_s(\shE)$
as the minimum of $\deg \, (\shQ)$, where
$\shQ$ is a quotient sheaf of rank $s$.
Since $\shG$ is ample, it follows that $\mindeg_s (\shG) >0$, hence
the estimate with 
$\deg \, (\shG)$ follows.
(This last estimate holds also for
$\bar{\mu}_{\rm max}$.)
\qed

\medskip
If the bundle has rank two, the following theorem of Hartshorne-Mumford
gives a satisfactory criterion for ampleness also in positive
characteristic.

\begin{theorem}
\label{amplecritranktwo}
Let $Y$ denote a smooth projective curve of genus $g$ defined
over an algebraically closed field of characteristic $p >0$.
Let $\shF$ denote a locally free sheaf of rank two on $Y$.
Suppose that $ \deg \, (\shF) > \frac{2}{p} (g-1)$ and
$\deg \, (\shL) >0$ for every invertible
quotient sheaf $ \shF \ra \shL \ra 0$.
Then $\shF$ is ample.
\end{theorem}

\proof
See \cite[Proposition 7.5 and Corollary 7.7]{hartshorneample}.
\qed

\begin{remark}
There exist more criteria
for ample and very ample vector bundles on curves,
see for example \cite{alzati}, \cite{campanaflenner}, \cite{ionescutoma}.
We omit them, since they
don't seem to have implications on tight closure problems.
\end{remark}

\section{Criteria for affineness}
\label{sectionaffine}

For ease of reference we fix the following situation.

\begin{situation}
\label{extensionsit}
Let $Y$ denote a smooth projective curve over an algebraically closed
field $K$ and let
$\shS$ denote a locally free sheaf on $Y$ and let
$\shG = \shS^\dual$ be its dual sheaf.
A cohomology class $c \in H^1(Y, \shS)= \Ext^1(\O_Y, \shS)$ yields
a short exact sequence
$$ 0 \lra \shS \lra \shS' \lra \O_Y \lra 0 \, .$$
The dual sequence $ 0 \ra \O_Y \ra \shG' \ra \shG \ra 0$
yields a projective subbundle $\PP(\shG) \subset \PP(\shG')$ of codimension
one.
\end{situation}

Such a situation arises in particular from a homogeneous
$R_+$-primary tight closure problem in a homogeneous
coordinate ring $R$ of $Y$.
This tight closure point of view leads to
the question whether the complement
$\PP(\shG') - \PP(\shG)$ is affine or not, in dependence of
$\shS$ and $c \in H^1(Y, \shS)$ (see section \ref{sectiongraded}).
Though we consider in this section mainly the case of characteristic $0$,
subsequent results in positive characteristic are discussed in section
\ref{sectionexclusion} and section \ref{sectionvanishing}.
The ampleness criterion \ref{amplekritnull}
yields at once the following affineness criterion.

\begin{corollary}
\label{ampleaffin}
Let the notation be as in {\rm \ref{extensionsit}} and suppose that
$\Char (K)=0$.
Suppose that $\shG$ is ample
\-- that is $\mu_{\rm min}(\shG) > 0$ or equivalently $\mu_{\rm max}(\shS) <0$
\--
and that $c \neq 0$.
Then $\PP(\shG')- \PP(\shG)$ is affine.
\end{corollary}

\proof
Since $\shG$ is ample and since $0 \ra \O_Y \ra \shG' \ra \shG \ra 0$
is a non splitting short exact sequence,
the sheaf $\shG'$ is also ample due to
\cite[Theorem 2.2]{giesekerample} (here we use characteristic $0$).
Hence $\PP(\shG) \subset \PP(\shG')$
is an ample divisor and its complement is affine.
\qed

\medskip
Even if the divisor $\PP(\shG) \subset \PP(\shG')$
is not ample, the open subset $\PP(\shG') - \PP(\shG)$ may be affine.
We need the following lemma to obtain more general sufficient criteria
for the affineness of the complement.

\begin{lemma}
\label{affinlemma}
Let $Y$ denote a scheme and let $\shS$ and $\shT$
be locally free sheaves on $Y$ and let
$q: \shS \ra \shT$ be a sheaf homomorphism.
Let $c \in H^1(Y,\shS)$ with corresponding extension
$0 \ra \shS \ra \shS' \ra \O_Y \ra 0$
and let $q(c) \in H^1(Y, \shT)$ be its image
with corresponding extension
$0 \ra \shT \ra \shT' \ra \O_Y \ra 0$.

If $\PP((\shT')^\dual) - \PP( \shT^\dual)$ is affine, then also
$\PP((\shS')^\dual)-\PP(\shS^\dual)$ is affine.

\end{lemma}
\proof
See \cite[Lemma 3.1]{brennertightelliptic}.
\qed

\begin{theorem}
\label{slopekritaffin}
Let the notation be as in {\rm \ref{extensionsit}} and suppose that
$\Char (K)=0$.
Suppose that there exists a semistable sheaf $\shT$
of negative slope, $\mu (\shT) < 0$,
and a sheaf morphism $q: \shS \ra \shT $
such that $q(c) \neq 0$ in $H^1(Y, \shT)$.
Then the complement $\PP(\shG') - \PP(\shG)$ is affine.
\end{theorem}
\proof
The sheaf $\shH= \shT^\dual$ is
semistable of positive degree, hence ample
due to \ref{amplekritsemistable},
therefore $\PP(\shH')-\PP(\shH)$ is affine due to \ref{ampleaffin}.
Therefore $\PP(\shG')-\PP(\shG)$ is also affine due to \ref{affinlemma}.
\qed

\medskip
The first candidates of semistable sheaves to look at 
are $\shS/ \shS_{s-1}$ (the semistable quotient of minimal slope)
and invertible sheaves.

\begin{corollary}
Let the notation be as in {\rm \ref{extensionsit}} and suppose that
$\Char (K)=0$.
Suppose that $\mu_{\rm min}(\shS) <0$
{\rm(}or equivalently $\mu_{\rm max}(\shG) >0${\rm )}
and that the image of
$c \in H^1(Y,\shS)$ in the semistable quotient $\shQ =\shS/ \shS_{s-1}$
is $\neq 0$.
Then $\PP(\shG') -\PP(\shG)$ is affine.
\end{corollary}

\proof
Since $\mu (\shQ)= \mu_{\min} (\shS) <0 $, the result
follows from \ref{slopekritaffin}.
\qed

\begin{corollary}
\label{quotientnegdegree}
Let the notation be as in {\rm \ref{extensionsit}} and
suppose that $\Char (K)=0$.
Suppose that $\shG$ is not normalized, that is
there exists an invertible
sheaf $\shL$ of negative degree such that
$H^0 (Y,\shG \otimes \shL) \neq 0$.
If $q: \shS \ra \shL$ is a corresponding non-trivial map
such that $0 \neq q(c) \in H^1(Y, \shL)$,
then $\PP(\shG') -\PP(\shG)$ is affine.
\end{corollary}

\proof
Note that $H^0(Y,\shG \otimes \shL) = \Hom_{\O_Y}(\shS, \shL)$,
thus this follows again from \ref{slopekritaffin}.
\qed

\begin{remark}
\label{ampleposremark}
In positive characteristic the assumption in \ref{ampleaffin}
that $c \neq 0$ is to weak to ensure the ampleness of $\shG'$.
We need the stronger condition that
$0 \neq \varphi^*(c) \in H^1(Y', \varphi^*(\shS))$
for every finite morphism $\varphi: Y' \ra Y$.
In the situation coming from a forcing problem in tight closure,
this condition means that $f_0$ does not belong to the plus closure.

However, if $\shG$ is invertible and of positive degree,
then the extension $\shG'$ corresponding to $0 \neq c \in H^1(Y, \shG^\dual)$
is ample also in positive characteristic $p$
under the condition that $p \geq 2 (g-1)$, where $g$ is the genus
of the curve $Y$.
This follows from the ampleness criterion
of Hartshorne-Mumford for bundles of rank two, see
\ref{amplecritranktwo}.

From this it follows also that Corollary \ref{quotientnegdegree}
holds also in positive characteristic $p \gg 0$.
For in this case the cohomology class $0 \neq q(c) \in H^1(Y, \shL)$
in the assumption of \ref{quotientnegdegree}
gives rise to an ample sheaf $(\shL')^\dual$ and then
the open subset $\PP( (\shL') ^\dual ) -\PP(\shL^\dual)$ is affine.
The affineness of $\PP(\shG') -\PP(\shG)$ follows then from \ref{affinlemma}.

If $\shG$ is indecomposable of rank $r \geq 2$,
then the degree of
$\shG$ must fulfill stronger conditions to ensure ampleness.
The result \cite[Theorem 25]{tango} suggests that the right inequality
may be $\deg \, (\shG) > r(r-1)(g-1)+ 2r (g-1)/p$
(or equivalently $\mu(\shG) > (g-1) (r+ p/2 -1)$).
Tango considers only the behavior of a cohomology class
$c \in H^1(Y,\shS)$ under the Frobenius, but along these lines
it should be possible to deduce also an ample criterion. 
\end{remark}

\begin{example}
Examples of semistable sheaves of positive degree
which are not ample (or which are not strongly semistable)
correspond to examples where the vanishing theorem for
tight closure does not hold in
small characteristic.
Look at the examples in \cite[\S3]{hartshorneamplecurve} and
\cite[Lemma 2.4]{giesekerample}) on the one hand and at
\cite[Example 3.11]{hunekesmithkodaira} on the other hand.
\end{example}

\section{Criteria for $\PP(\shG') -\PP(\shG)$ not to be affine}
\label{sectionnonaffine}

We look again at the situation of \ref{extensionsit},
but now we look for criteria for $\PP(\shG') -\PP(\shG)$ to be not affine.
Since affine subsets contain no projective curves,
we get the following easy criterion.

\begin{proposition}
Let the notation be as in {\rm \ref{extensionsit}}.
Suppose that there exists a finite morphism
$\varphi: Z \ra Y$ such that
$\varphi^*(c)= 0$. Then $\PP(\shG') - \PP(\shG)$ is not affine.
\end{proposition}
\proof
The condition means that the pull back of
$ 0 \ra \O_Y \ra \shG' \ra \shG \ra 0$ splits on $Z$.
The splitting surjection $  \varphi^*(\shG') \ra \O_Z \ra 0$
yields a section $Z \ra \PP(\varphi^*(\shG'))$ disjoined to
$\PP(\varphi^*(\shG))$ and hence there exists a projective curve
inside $\PP(\shG') - \PP(\shG)$.
\qed

\begin{remark}
\label{bigaffinremark}
If the complement of an effective Weil divisor $D$ on a normal projective
variety $X$ of dimension $d$ is affine, then
$\Gamma(X-D, \O_X)$ is a finitely generated $K$-algebra of dimension $d$.
Thus the Iitaka-dimension of $D$ is maximal, and $D$ is called big,
see \cite[Definitions 2.1.3 and 2.2.1]{lazarsfeld}
This means that there exists a number $c >0$ such that
$h^0(X, \O_X(kD)) \geq ck^d$ for $k \gg 0$.
Hence from conditions which imply that
$S^k(\shG)= \Gamma (\PP(\shG), \O_{\PP(\shG)}(1))$ has few sections
we may derive that $\PP(\shG') - \PP(\shG)$ is not affine.
\end{remark}

\begin{lemma}
\label{muminpseudo}
Let $Y$ denote a smooth projective curve over an algebraically closed
field and let $\shT$ denote a locally free sheaf on $Y$.
Suppose that $\bar{\mu}_{\min} (\shT) \geq 0$.
Then $\mu_{\min} ( \Gamma^k (\shT)) \geq 0$ for $k \geq 0$, where
$\Gamma^k( \shT) = (S^k(\shT^\dual))^\dual $.
\end{lemma}
\proof
Note that $S^k(\shT) \cong (S^k(\shT^\dual))^\dual$
is only true in characteristic zero.
Assume that there exists a locally free quotient sheaf
$(S^k (\shT^\dual))^\dual \ra \shQ \ra 0$ of negative degree.
We find a finite morphism
$\varphi: Y' \ra Y$ such that we may write
$\varphi^*(\shQ) = \shL^k \otimes \shN $, where
$\shL$ is an invertible sheaf with $\deg \, (\shL) < 0$ and also
$ \deg \, (\shN) <0$.
Due to our assumption we may assume that this is already true on $Y$.
We tensor by $\shL^{-k}$ and get
$(S^k (\shT^\dual))^\dual \otimes \shL^{-k} \ra \shN \ra 0$.
But
$$(S^k (\shT^\dual))^\dual \otimes \shL^{-k}
=(S^k(\shT^\dual) \otimes \shL^k)^\dual
=(S^k(\shT^\dual \otimes \shL)   )^\dual
=(S^k( (\shT \otimes \shL^\dual)^\dual)^\dual \,.$$
Now $\deg \, (\shL^\dual) >0$, hence
$\bar{\mu}_{\min} (\shT \otimes \shL^\dual) >0$ and therefore
$\shT \otimes \shL^\dual$ is ample due to \ref{amplekrit}.
Then due to \cite[Theorem 6.6 and Proposition 7.3]{hartshorneample}
it follows that $(S^k( (\shT \otimes \shL^\dual)^\dual)^\dual$ is ample,
buts its quotient sheaf $\shN$ is not, since $\deg \, (\shN) <0$,
which gives a contradiction.
\qed

\begin{theorem}
\label{slopemaxkrit}
Let the notation be as in {\rm \ref{extensionsit}}.
Suppose that
$\bar{\mu}_{\rm max} (\shG ) \leq 0$. 
Then $\PP(\shG')- \PP(\shG)$ is not affine.
\end{theorem}

\proof
From the sequence $0 \ra \O_Y \ra \shG' \ra \shG \ra 0$ it follows that
$\bar{\mu}_{\max} (\shG') \leq 0$ holds as well, since the image of
a mapping $\shS \ra \shG'$,
where $\shS$ is a semistable sheaf of positive slope,
must lie inside the kernel of $\shG' \ra \shG$.
Applying Lemma \ref{muminpseudo} to the dual of $\shG'$
it follows that $\mu_{\max} (S^k(\shG')) \leq 0$.

If moreover $ \mu_{\rm max} (S^k (\shG')) < 0$,
then $\Gamma(Y, S^k(\shG'))=0$.
So we may suppose that $ \mu_{\rm max} (S^k (\shG')) = 0$.
Then there exists the maximal destabilizing sheaf
$0 \ra \shF_k \ra S^k (\shG') \ra S^k (\shG')/ \shF_k \ra 0$
such that
$\mu (\shF_k)=0$ and $\mu_{\rm max}(S^k(\shG')/\shF_k) <0$,
so that this sheaf has again no global sections $\neq 0$.

We claim that for a semistable locally free sheaf $\shF$ of degree zero and
rank $s$ the dimension of the global sections is at most $s$:
This is true for invertible sheaves so we do induction on the rank.
If $\shF$ has a global section $\neq 0$, then $\O_Y \subseteq \shF$
is a subsheaf.
We consider the saturation $\O_Y \subseteq \shM \subseteq \shF$
so that the cokernel $ \shF /\shM$ is torsion free, hence locally free on the
curve (see \cite[1.1]{huybrechtslehn}).
$\shM$ has degree $0$ (since $\O_Y \subseteq \shM$ and
since $\shF$ is semistable) and therefore we may apply the induction hypothesis
to the cokernel.

This gives the estimates (set $r= \rk (\shG)$ and $r+1= \rk (\shG')$)
$$
h^0 (S^k (\shG')) \leq h^0(\shF_k) + h^0(S^k (\shG')/\shF_k)
\leq {\rm rk} (\shF_k) \leq {\rm rk} (S^k (\shG'))= 
\binom{k+r}{r} \, . $$
This is a polynomial with leading coefficient
$\frac{1}{r !} k^{r}$, thus it is bounded above
by $\leq c k^{r}$ for some $c >0$.
Since the dimension of $\PP(\shG')$ is $r+1$ and
since $\pi_* \O_{\PP(\shG')}(k) =S^k(\shG')$ it follows that
$\O_{\PP(\shG')}(1)$ is not big and that $\PP(\shG')-\PP(\shG)$
is not affine by remark \ref{bigaffinremark}.
\qed

\begin{corollary}
Let the notation be as in {\rm \ref{extensionsit}}.
Suppose that $\shG$ is strongly semistable and $\mu(\shG) \leq 0$.
Then $\PP(\shG')- \PP(\shG)$ is not affine.
\end{corollary}

\proof
Since $\shG$ is strongly semistable we have
$\bar{\mu}_{\max} (\shG) = \mu(\shG) \leq 0$, hence
the result follows from \ref{slopemaxkrit}.
\qed

\begin{remark}
\label{remarkelliptic}
For an elliptic curve $Y$ one can improve the Theorems
\ref{slopekritaffin} and \ref{slopemaxkrit}
to a complete characterization of affine subsets (in any characteristic).
Suppose the notation of \ref{extensionsit}
and let $\shS= \shS_1 \oplus \ldots \oplus \shS_s$ be the decomposition
into indecomposable sheaves, $c=(c_j)$, $c_j \in H^1(Y,\shS_j)$.
An indecomposable sheaf on an elliptic curve is strongly semistable,
hence we have in this situation a much better control on the different
notions of slopes.
Theorem 3.2 of \cite{brennertightelliptic} states that
$\PP(\shG') - \PP(\shG)$ is affine if and only if
there exists $j$ such that $\deg \, (\shS_j) <0$ and $c_j \neq 0$.
In positive characteristic, the same numerical condition is equivalent
to the non-existence of projective curves inside $\PP(\shG') - \PP(\shG)$,
see \cite[Theorem 3.3]{brennertightelliptic}.
\end{remark}

\section{Projective bundles corresponding to tight closure problems}
\label{sectiongraded}

In this section we briefly recall some results from
\cite{brennertightproj} on how graded tight closure problems
in a graded ring $R$ translate to
problems about projective bundles and subbundles over $\Proj\, R$.
Let $K$ denote an algebraically closed field and let $R$ be a
standard $\NN$-graded $K$-algebra, that is $R_0=K$ and $R$ is generated
by finitely many elements of degree one. Set $Y= \Proj \, R$.
Let $f_i$ be homogeneous $R_+$-primary elements of $R$ of degree $d_i$,
that is the $D_+(f_i)$ cover $Y$. Then for $m \in \ZZ$
we have the exact sequence of locally free sheaves on $Y$,
$$ 0 \lra \shR(m) \lra \O_Y(m-d_1) \oplus \ldots \oplus \O_Y(m-d_n)
\stackrel{ \sum f_i}{\lra} \O_Y( m) \lra 0 \, ,$$
and dually
$$0 \lra \O_Y(-m) \lra \O_Y(d_1-m) \oplus \ldots \oplus \O_Y(d_n-m)
\lra \shF(-m) \lra 0 \, .$$
$\shR(m)$ is the {\em sheaf of relations} of total degree $m$ corresponding
to the generators $f_1, \ldots, f_n$. Its rank is $n-1$.
From this {\em presenting sequence} it follows
that $\Det ( \shF(-m)) = \O_Y(d_1+ \ldots +d_n - (n-1)m)$
and $\deg \, (\shF(-m))= (d_1+ \ldots +d_n - (n-1)m) \deg \, (\O_Y(1))$.
A geometric realization of the vector bundle
$ V(-m) = \Spec \, \oplus S^k (\shF(-m))$
is given by the open subset
$$ V(-m) = D_+(R_+)
\subset \Proj\, R[T_1,\ldots,T_n]/ (\sum_{i=1}^n f_iT_i) \, ,$$
where $\deg \, (T_i) = m-d_i$.

\smallskip
Now suppose that $f_0$ is another homogeneous element in $R$ of degree
$m$. The elements $f_1, \ldots ,f_n,f_0$ define again a sheaf of relations
$\shR'(m)$ together with a short exact sequence
$$0 \lra \shR(m) \lra \shR'(m) \lra \O_Y \lra 0 \, ,$$
which we call the {\em forcing sequence}.
The corresponding cohomology class $c \in H^1(Y, \shR(m))$
(the {\em forcing class})
is also given by the connecting homomorphism
$R_m \ra \Gamma(Y, \O_Y(m)) \ra H^1(Y, \shR(m))$.
The dual sequence
$$ 0 \lra \O_Y \lra \shF'(-m) \lra \shF(-m) \lra 0 $$
yields a closed subbundle $\PP(\shF(-m)) \subset \PP(\shF'(-m))$
(we will skip the number $m$ in this expression and write
$\PP(\shF)$).

The basic fact is that
the complement $\PP(\shF') -\PP(\shF)$
is isomorphic to the $\Proj$
of the so called forcing algebra
$R[T_1, \ldots, T_n]/(f_1T_1 + \ldots +f_nT_n +f_0)$
(suitable graded).
See \cite{hochstersolid} for forcing algebras and how
the tight closure of an ideal is expressed in terms of them.
The containment of a homogeneous element in the ideal,
in the tight closure and in the plus closure of the ideal
is expressed in terms of these projective bundles in the following way.
Note that in characteristic zero we work with the notion of solid closure.

\begin{lemma}
\label{trivialtwo}
In the described situation the following
are equivalent.

\renewcommand{\labelenumi}{(\roman{enumi})}
\begin{enumerate}

\item 
$f_0 \in (f_1, \ldots ,f_n)$.

\item
There exists a section $Y \ra \PP(\shF')$ disjoined to
$\PP(\shF) \subset \PP(\shF')$.

\item
The forcing sequence
$0 \ra \shR(m) \ra \shR'(m) \ra \O_Y \ra 0 $ splits.

\item
The corresponding cohomological class in $H^1(Y,\shR(m))$ vanishes.
\end{enumerate}
\end{lemma}

\proof
See \cite[Lemma 3.7]{brennertightproj}.
\qed

\begin{proposition}
\label{geokrittight}
Let $R$ be a normal standard-graded $K$-algebra of dimension $2$,
let $f_1, \ldots ,f_n \in R$ be $R_+$-primary homogeneous elements
and let $f_0$ be another homogeneous element.
Then $f_0 \in (f_1, \ldots ,f_n)^*$ if and only if
$\PP( \shF') - \PP(\shF)$ is not affine.

Furthermore, if the characteristic of $K$ is positive,
the following are equivalent.
\renewcommand{\labelenumi}{(\roman{enumi})}
\begin{enumerate}

\item
$f_0 \in (f_1, \ldots ,f_n)^{+{\rm gr}}$, that is there exists
a finite graded extension $R \subseteq R'$ such that
$f_0 \in (f_1, \ldots ,f_n)R'$.

\item
There exists a smooth projective curve $Z$
and a finite surjective morphism $g: Z \ra Y$ such that
the pull back $g^*\PP(\shF')$
has a section not meeting $g^*\PP(\shF)$.

\item
There exists a curve $Z \subset \PP(\shF')$
which does not intersect $\PP(\shF)$.

\end{enumerate}
\end{proposition}
\proof
See \cite[Lemmata 3.9 and 3.10]{brennertightproj}.
\qed

\section{Applications to tight closure: inclusion bounds}
\label{sectioninclusion}

We shall now apply the results of the previous sections to tight
closure problems. We fix the following situation.

\begin{situation}
\label{forcingsit}
Let $K$ denote a field with algebraic closure $\bar{K}$ and let
$R$ denote a two-dimensional standard-graded
$K$-algebra such that $R_{\bar{K}}$ is a normal domain.
Let $Y= \Proj\, R$ denote the cor\-responding
smooth projective curve over $\bar{K}$, let $g$ denote its genus
and let $\delta = \deg \, (Y) = \deg \, (\O_Y(1))$ denote the degree
of the very ample invertible sheaf $\O_Y(1)$ on $Y$.
Let $f_1, \ldots ,f_n \in R$ be homogeneous $R_+$-primary elements
of degree $d_i$.
Let $\shR(m)$ be the sheaf of relations
of total degree $m$ on $Y$ and let $\shF(-m)$ be its dual sheaf.

Let $f_0$ denote another homogeneous element of degree $m$,
let $\shR'(m)$ be the sheaf of relations for the elements
$f_1, \ldots, f_n,f_0$ of total degree $m$ and let
$$ 0 \lra \shR(m) \lra \shR'(m) \lra \O_Y \lra 0 $$
be the corresponding extension and
let $c \in H^1(Y, \shR(m))$ be the corresponding forcing class
defined by $f_0$.
The corresponding surjection
$\shF'(-m) \ra \shF(-m) \ra 0$ yields the embedding
$\PP(\shF(-m)) \subset \PP(\shF'(-m))$.
\end{situation}

\begin{definition}
Let the notation be as in \ref{forcingsit}.
Then we set
$$ \mu_{\rm max}(f_1, \ldots ,f_n) := \mu_{\rm max} (\shF(0)) \, $$
and also for $\mu$, $\mu_{\min}$, $\bar{\mu}_{\max}$ and $\bar{\mu}_{\min}$.
\end{definition}

\begin{remark}
Note that we consider for ideal generators $f_1, \ldots ,f_n \in R$
always the slope of the relation bundle after replacing
$K$ by $\bar{K}$. Since changing the base field does not affect
the affineness of open subsets, it does not affect solid closure.

The slope of $\shF(0)$ is \-- due to the presenting sequence for $\shF(0)$ \--
given by
$$\mu(\shF(0)) = \frac{d_1+ \ldots +d_n}{n-1} \delta  \, ,$$
therefore we get the estimates
$$ \mu_{\rm min}(f_1, \ldots ,f_n) \leq
\frac{d_1+ \ldots +d_n}{n-1} \delta \leq 
\mu_{\rm max}(f_1, \ldots ,f_n) \, .$$
Equality holds if and only if $\shF(0)$ is semistable.
Furthermore we have
$$\mu(\shF(-m))=\mu (\shF(0) \otimes \O_Y(-m)) =  \mu(\shF(0)) -m \delta 
= ( \frac{d_1+ \ldots +d_n}{n-1} -m )\delta  \, ,$$
and the same rule holds for $\mu_{\rm max}$ etc.
If $f_0$ is another homogeneous element of degree $m$, then the number
$\deg \, (\shF(-m))= \deg \, (\shF'(-m)) = (d_1+ \ldots +d_n -(n-1)m) \delta$
is also the top self intersection number of the forcing subbundle
$\PP(\shF) \subseteq \PP(\shF')$.
\end{remark}

From the conditions in section \ref{sectionnonaffine}
we derive the following
numerical condition that elements of sufficiently
high degree must belong to the tight closure.

\begin{theorem}
\label{maxin}
Let the notation be as in {\rm \ref{forcingsit}}. If
$ \deg \,(f_0) \geq  \frac{1}{\delta}  \bar{\mu}_{\rm max} (f_1, \ldots ,f_n)$,
then $f_0 \in (f_1, \ldots,f_n)^*$.
\end{theorem}
\proof
Let $m = \deg \, (f_0)$.
The condition means that
$$ \bar{\mu}_{\rm max} (\shF (-m))
=\bar{\mu}_{\rm max} (f_1, \ldots, f_n) -m \delta \leq 0 \, .$$
Hence the result follows from \ref{slopemaxkrit}
and \ref{geokrittight}.
\qed

\begin{remark}
This numerical criterion generalizes
the corresponding statement for cones over elliptic curves
proved in \cite[Corollary 4.9]{brennertightelliptic}.
On an elliptic curve we have the equality
$\bar{\mu}_{\max} (\shF(0)) = \mu_{\max} (\shF(0)) = \max_j (\mu(\shF_j))$,
where $\shF(0)= \oplus_j \shF_j$ is the decomposition into
indecomposable sheaves.
\end{remark}

To obtain criteria for tight closure membership
we need bounds from above for $\bar{\mu}_{\max} (f_1, \ldots ,f_n)$.
The next proposition gives a general bound for $\bar{\mu}_{\rm max}$.
We will give a much better bound in section \ref{examples} for the case $n=3$
under the condition that the sheaf of relations is indecomposable.

\begin{proposition}
\label{mumaxabschaetzung2}
Let the notation be as in {\rm \ref{forcingsit}}.
Suppose that the degrees are ordered
$1 \leq d_1 \leq d_2 \leq \ldots \leq d_n$.
Let $\shE= \oplus \O_Y(d_i)$ and let
$0 \ra \O_Y \ra \shE \ra \shF(0) \ra 0$ be the presenting sequence
for $\shF(0)$.
Then we have the estimate
$$ \mu_{\rm max}(f_1, \ldots ,f_n)
\leq
\delta \cdot (d_{n-1}+d_n) \, .$$
The same is true for $\bar{\mu}_{\max}$.
\end{proposition}

\proof
Set $\shF= \shF(0)$. Corollary
\ref{mumaxabschaetzung} together with the inequality
$\mindeg_s(\shF) \geq \mindeg_s(\shE)$
yields
$$
\mu_{\rm max} (\shF)
\leq
\max_{s=0, \ldots ,n-2} \frac{ \deg \, (\shF) -\mindeg_s\, (\shF) }{n-1-s}
\leq
\max_{s=0, \ldots ,n-2}
\frac{ \deg \, (\shE) -\mindeg_s\, (\shE) }{n-1-s} \, .
$$
We claim that
$\mindeg_s \, (\shE) = \delta(d_1 + \ldots +d_s)$.
Since $\O_Y(d_1) \oplus \ldots \oplus \O_Y(d_s)$
is a quotient sheaf of rank $s$, the estimate $\leq$
is clear.
For the other estimate we consider first the case
$s=1$, so suppose that $\shQ$ is an invertible sheaf.
If $\shQ$ is a quotient of $\shE$, then
$\Hom( \O_Y(d_i), \shQ)=H^0(Y, \O_Y(-d_i) \otimes \shQ) \neq 0$ for
at least one $i$. Therefore
$\deg \, (\shQ) \geq \delta {\rm min}_i d_i = \delta d_1$.

Now suppose that $\shQ$ is a locally free quotient of $\shE$ of rank $s$.
Then we have a surjection
$$\oplus_{i_1 < \ldots < i_s}\,
\O_Y(d_{i_1}) \otimes \ldots \otimes \O_Y(d_{i_s})
\cong
\bigwedge^s \shE
\lra \bigwedge^s \shQ= \det \shQ \, .
$$
Due to the case $s=1$ we know
$\deg \, (\shQ) \geq \delta (d_1 + \ldots + d_s)$, which proves the claim.

Thus we have the estimate
$$
\mu_{\rm max} (\shF)
\leq
\max_{s=0, \ldots ,n-2}
\frac{ \deg \, (\shE) -\mindeg_s\, (\shE) }{n-1-s}
=\max_{s=0, \ldots ,n-2}
\delta \frac{d_{s+1}+ \ldots + d_n}{n-1-s}
 \, .
$$
Here the term for $s=n-2$, which is $\delta ( {d_{n-1} +d_n})$,
is maximal.

If $\varphi: Y' \ra Y$ is a finite morphism,
then the situation is preserved under the pull-back
(even if $\varphi^* \O_Y(1)$ is not very ample anymore).
Then
${\mu}_{\max}(\varphi^*(\shF))
\leq (\delta \cdot \deg \, (\varphi)) (d_{n-1} + d_n) $
and hence the inequality holds also for
$\bar{\mu}_{\max}$.
\qed

\begin{corollary}
\label{inclusionbound}
Let the notation be as in {\rm \ref{forcingsit}}.
Suppose that the degrees are ordered
$1 \leq d_1 \leq d_2 \leq \ldots \leq d_n$.
Then for $m \geq d_{n-1} + d_n$ we have the inclusion
$$ R_{m} \subseteq (f_1, \ldots ,f_n)^*  \, .$$
\end{corollary}
\proof
This follows from \ref{mumaxabschaetzung2} and \ref{maxin}.
\qed

\begin{remark}
This corollary is for the two-dimensional case
a somewhat better estimate than
the estimate $2 \max_i (d_i)$
(and $d_1+ \ldots + d_n$) found by
K. Smith, see \cite[Proposition 3.1 and Proposition 3.3]{smithgraded}.
Smith´s estimate for $R_+$-primary ideals
is however true in every dimension
and holds also for the plus closure
in positive characteristic.
\end{remark}

\begin{remark}
An estimate from below for $\mu_{\max}$ is
$\mu_{\rm max} (\shF) \geq \delta \cdot {\rm max}_i (d_i)$.
To see this consider again the sequence
$0 \ra \O_Y \ra \oplus_i\,  \O_Y(d_i) \ra \shF \ra 0$.
If $\O_Y(d_i) \ra \shF $ is the zero map for one $i$, then
$\O_Y \cong \O_Y(d_i)$, $d_i=0$ and $f_i$ is a unit and the statement is
clear from $\shF \cong \oplus_{j \neq i} \O_Y(d_j)$.
Otherwise $\O_Y(d_i) \ra \shF $ is not the zero map, and then
$\mu_{\rm max} (\shF) \geq \mu(\O_Y(d_i)) = \delta d_i$ for all $i$.
\end{remark}

\section{Applications to tight closure: exclusion bounds}
\label{sectionexclusion}

We are now looking for degree bounds $a$ such that if $\deg (f_0) <a$,
then $f_0 \in (f_1,...,f_n)^*$ if and only if $f_0 \in (f_1, \ldots, f_n)$.
So below the degree bound an element $f_0$ is excluded from the tight closure
with the exception that it belongs to the ideal itself.

The theorems in this and the next section hold either in characteristic zero
or in positive characteristic $p$ under the condition that $p  \gg 0$.
To make sense of this statement we have to
suppose that everything is given relatively to a base scheme
such that the generic fiber has characteristic zero and the special fibers
have positive characteristic.
For this we fix the following situation, see
also \cite[Definition 3.3]{hunekesmithkodaira}
and in particular the appendix of Hochster in \cite{hunekeapplication}
for this setting.

\begin{situation}
\label{relsituation}
Let $D$ denote a finitely generated normal $\ZZ$-domain of dimension one.
Let $S$ denote a standard-graded flat $D$-algebra
such that for all $\fop \in \Spec \, D$ the algebras
$S_{\kappa(\fop)} = S \otimes_D \kappa(\fop) $
are two-dimensional geometrically normal standard-graded
$\kappa(\fop)$-algebras (so that $S_{\bar{\kappa}(\fop)}$ are normal domains).
For $\fop =0$ this is an algebra over the quotient field $Q(D)$ of
characteristic zero, and for a maximal ideal $\fop$
the algebra $S_{\kappa(\fop)}$ is an algebra over the finite field
$\kappa(\fop)=D/ \fop$ of positive characteristic.

We suppose that we have $S_+$-primary homogeneous elements
$f_1, \ldots, f_n \in S$ of degree $d_i$ and another
homogeneous element $f_0$.
Let $B$ denote the forcing algebra over $S$ for this data and let
$U= D(S_+) \subseteq \Spec \, B$.
These elements yield homogeneous forcing data for every $S_{\kappa(\fop)}$
and $B_{\kappa(\fop)}=B \otimes_D \kappa(\fop)$
is the corresponding forcing algebra.
For every prime ideal $\fop \in \Spec\, D$ the affineness of
$U_\fop =U_{\kappa(\fop)}= U \cap \Spec \, B_{\kappa(\fop)}$
is equivalent to $f_0 \not\in (f_1, \ldots ,f_n)^*$ in $S_{\kappa(\fop)}$.

We denote by $Y= \Proj\, S$ the smooth projective (relative) curve
over $\Spec\, D$
and by $\delta$ the common degree of the curves $Y_\fop$, $\fop \in \Spec\, D$.
We denote by $\shR(m)$ the sheaf of relations on $Y$ and by $\shF(-m)$
its dual sheaf and we denote the restrictions to $Y_\fop$
by $\shR_\fop(m)$ and $\shF_\fop(-m)$.
The different notions of slopes and of semistability refer always
to $\shF_{\bar{\kappa}(\fop)}$ on $Y_{\bar{\kappa}(\fop)}$.
The element $f_0$ yields an extension
$0 \ra \shR(m) \ra \shR'(m) \ra \O_Y \ra 0$ and a subbundle
$\PP(\shF) \subset \PP(\shF')$, which induces the projective subbundle
$\PP(\shF_{\fop}) \subset \PP(\shF'_{\fop})$ on every curve
for every $\fop \in \Spec\, D$.
\end{situation}

\begin{remark}
\label{affineopen}
We will apply several times the following conclusion:
let $D \subseteq B$ denote Noetherian domains and let
$U=D(\foa) \subseteq \Spec\, B$ denote an open subset,
$\foa=(a_1, \ldots ,a_k)$.
Suppose that $U_\eta = U \cap \Spec \, (B \otimes_D \kappa(\eta))$
is affine, where $\eta$ denotes the generic point of $\Spec \, D$.
This means that there exist rational functions
$q_j \in \Gamma(U_\eta , \O_{\eta})$ such that $\sum q_ja_j =1$.
We find a common denominator $0 \neq g \in D$ such that
these functions $q_j$ are defined on $U \cap D(g) \subseteq \Spec\, B$,
hence also $U \cap D(g)$ is affine.
This means that after shrinking $D$ (i. e. replacing $\Spec D$ by $\Spec D_g$)
we may assume that $U$ itself is affine.
Hence for every $P \in \Spec \, D$ the fibers $U_{\kappa(P)}$ are affine.

Suppose in the situation \ref{relsituation}
that $f_0 \not\in (f_1, \ldots ,f_n)^*$ holds over the generic point
$\eta \in \Spec \, D$. This means that the open subset
$U_\eta$ is affine.
Then after shrinking $D$ we may assume that $U$ is affine,
hence that every fiber $ U_{\kappa(\fop)}$ is affine.
This means that $f_0 \not\in (f_1, \ldots ,f_n)^*$ holds in $S_{\kappa(\fop)}$
for all $\fop \in \Spec \, D$ (or for almost all $\fop \in \Spec \, D$
for the old $D$). In this case we say briefly that
$f_0 \not\in (f_1, \ldots ,f_n)^*$ holds for $p \gg 0$.
\end{remark}

\begin{theorem}
\label{minex}
Let the notation be as in {\rm \ref{forcingsit}} and in
{\rm \ref{relsituation}}.
Suppose that the characteristic of $K$ is $0$ or $p  \gg  0$.
If $\deg \, (f_0) < \frac{1}{\delta} \mu_{\rm min} (f_1, \ldots ,f_n)$
{\rm(}$=\frac{1}{\delta} \mu_{\rm min}(\shF_\eta (0))${\rm )},
then $f_0 \in (f_1, \ldots, f_n)^*$
if and only if $f_0 \in (f_1, \ldots ,f_n)$.
\end{theorem}
\proof
Let $m =\deg \, (f_0)$. Suppose first that the characteristic is zero.
We may assume that $K$ is algebraically closed.
The condition means that
$ \mu_{\rm min}(\shF (-m) ) 
= \mu_{\rm min} (f_1, \ldots, f_n) -m \delta > 0$,
hence $\shF(-m)$ is ample due to \ref{amplekritnull}.
Suppose that $f_0 \not\in (f_1, \ldots ,f_n)$.
This means by Lemma \ref{trivialtwo} that the
corresponding forcing class is $c \neq 0$.
Hence $\PP(\shF') -\PP(\shF)$ is affine due to \ref{ampleaffin} and
$f_0 \not\in (f_1, \ldots ,f_n)^*$ due to \ref{geokrittight}.

Now suppose the relative situation \ref{relsituation}.
Note that the slope condition is imposed on the
generic fiber.
We have to show that $f_0 \not\in (f_1, \ldots ,f_n)$
implies $f_0 \not\in (f_1, \ldots ,f_n)^*$ for almost
all $\fop \in \Spec\, D$.
From $f_0 \not\in (f_1, \ldots ,f_n)$ in $S_{\kappa(\fop)}$
it follows $f_0 \not\in (f_1, \ldots ,f_n)$ in $S$
and by shrinking $D$ we may assume that
$f_0 \not\in (f_1, \ldots ,f_n)$ in $S_{Q(D)}$.
From the case of characteristic zero
we know that $U_\eta $ is affine
and the result follows from remark \ref{affineopen}.
\qed

\medskip
From the bounds proved in Theorem \ref{maxin} and in
Theorem \ref{minex} it is easy to derive the following result
of Huneke and Smith (see \cite[Theorem 5.11]{hunekesmithkodaira}).

\begin{corollary}
\label{projdim}
Let the notation be as in {\rm \ref{forcingsit}}
{\rm(}or {\rm \ref{relsituation} )}.
Suppose that the characteristic of $K$ is $0$
or $p  \gg  0$.
Suppose that the projective dimension of $R/(f_1, \ldots, f_n)$
is $2$ or equivalently that the sheaf of relations
$\shR(0) \cong \O_Y(-a_1) \oplus \ldots \oplus \O_Y(-a_r)$.
Let $a= \max \{ a_i\}$ and $b= \min \{a_i \}$.
Then
$$ R_{\geq a} \subseteq (f_1, \ldots ,f_n)^* 
\mbox{ and }
(f_1, \ldots ,f_n)^*  \subseteq (f_1, \ldots , f_n) + R_{ \geq b}\, .$$
\end{corollary}
\proof
Whenever $\shF(0)$ is a direct sum of invertible sheaves $\shL_j$ we have
$\bar{\mu}_{\rm max} (\shF(0)) = \max_j \deg \, (\shL_j)$
and $\bar{\mu}_{\min} (\shF(0)) = \min_j \deg \, (\shL_j)$.
So for $\shF(0)= \O_Y(a_1) \oplus \ldots \oplus \O_Y(a_r)$ we find
that $\bar{\mu}_{\max}(\shF(0)) =a \delta$ and
$\bar{\mu}_{\min} (\shF(0))= b \delta $,
so the result follows from \ref{maxin} and \ref{minex}.
\qed

\medskip
If the sheaf of relations on the projective curve splits into
invertible sheaves as in the previous Corollary
\ref{projdim} it is easy to give a numerical criterion for tight closure.

\begin{theorem}
\label{splitting}
Let the notation be as in {\rm \ref{forcingsit}}
{\rm(}or {\rm \ref{relsituation})}.
Suppose that we have a decomposition
$\shR(0)= \shL_1 \oplus \ldots \oplus \shL_s$,
where $\shL_j$ are invertible sheaves of degree $\mu_j$.
Let $m= \deg\, (f_0)$ and let $c \in H^1(Y, \shR(m))$
be the forcing class with components
$c_j \in H^1(Y ,\shL_j \otimes \O_Y(m))$.
Suppose that the characteristic is zero or $p \gg 0$.
Then $f_0 \in (f_1,\ldots ,f_n)^*$
if and only if $ \,\deg \, (\shL_j) +m \delta \geq 0 $ or $c_j=0$ holds
for all $1 \leq j \leq s$.
\end{theorem}
\proof
Suppose first that $f_0 \in (f_1, \ldots,f_n)^*$ holds and consider to the
contrary that there exists $j$ such that
$\deg \, (\shL_j) +m\delta < 0$ and $c_j \neq 0$.
Then Corollary \ref{quotientnegdegree} together with remark \ref{ampleposremark}
yields the contradiction $f_0 \not\in (f_1, \ldots,f_n)^*$.

For the other direction we consider the direct summand
$$\shS = \oplus_{\deg \, (\shL_j (m)) \geq 0}\, \shL_j(m)
\subseteq \shR(m) \, .$$
Then all non-zero components of the forcing class $c$
belong to $\shS$, so this class comes from and goes to a cohomology class
in $H^1(Y, \shS)$. Hence we know by Lemma \ref{affinlemma}
that the affineness of $\PP((\shS')^\dual)-\PP(\shS^\dual)$
is equivalent to the affineness of $\PP(\shF') - \PP( \shF)$.
Since the degree of every invertible summand sheaf of $\shS$ is nonnegative,
we know that $\bar{\mu}_{\rm min} (\shS) \geq 0$
and then $\bar{\mu}_{\rm max} (\shS^\dual) \leq 0$, therefore
$\PP((\shS')^\dual)-\PP(\shS^\dual)$ is not affine due to
Theorem \ref{slopemaxkrit}.
\qed

\begin{remark}
The situation of Theorem \ref{splitting} holds for every primary
homogeneous ideal in $K[x,y]$ (due to the splitting theorem of Grothendieck,
see \cite[Theorem 2.1.1]{okonekvector}),
but for a polynomial ring the computation of tight closure does not make
much problems, so \ref{splitting} gives a help which we do not need in
this case. However, the splitting situation holds also if
$I \subseteq R$ is the extended ideal $I=JR$ of an ideal
$J \subseteq K[x,y] \subset R$, and in this case it is also useful for
computations, see example \ref{indepexample} below.
There is also a version of Theorem \ref{splitting}
if the sheaf of relations splits into a direct sum of (strongly)
semistable sheaves.
\end{remark}

\begin{remark}
Let the notation be as in \ref{forcingsit}. From the sequence
$$ 0 \lra \O_Y \lra \O_Y(d_1) \oplus \ldots \oplus \O_Y(d_n) \lra
\shF \lra 0 \, $$ we also get the estimate
$$\bar{\mu}_{\rm min} (\shF) \geq \bar{\mu}_{\rm min} (\oplus_i \O_Y(d_i))
={\rm min}_i\, \mu(\O_Y(d_i))
= \delta {\rm min}_i (d_i) \, .$$
So if $f_0 \neq 0$ and
$ \deg \, (f_0) < {\rm min}_i (d_i)$, then $f_0 \not\in (f_1, \ldots , f_n)$
and due to \ref{minex} also $f_0 \not\in (f_1, \ldots , f_n)^*$.
\end{remark}

The minimal slope in the generic point gives also a bound
for the minimal slope in positive characteristic for $p \gg 0$.
The same is true for the maximal slope.

\begin{proposition}
\label{compare}
Suppose the relative situation \ref{relsituation},
let $Y= \Proj \, S \ra \Spec \, D$ denote the smooth projective
curve of relative dimension one and suppose
that the generic curve $Y_\eta=Y_{Q(D)}$ has at least one $Q(D)$-rational point.
Let $\mu_{\rm min} (\shF_\eta (-m))$
denote the minimal slope and let $\mu_{\rm max}(\shF_\eta (-m)) $
denote the maximal slope on $Y_{\bar{\kappa}(\eta)}$.
Then for $p \gg 0$ we have the bounds
$$ \bar{\mu}_{\rm min} (\shF(-m))
> \lceil \mu_{\rm min} (\shF_\eta(-m))  \rceil -1 \, $$
and
$$ \bar{\mu}_{\rm max} (\shF(-m))
< \lfloor \mu_{\rm max}(\shF_\eta (-m))  \rfloor +1 \, . $$
\end{proposition}
\proof
First we may assume by shrinking $D$ that there exists a section
for the relative curve $Y$. Hence there exists an invertible sheaf $\shM$
on $Y$ such that the degree of $\shM$
on every fiber $Y_\fop$ is one.
Let $\shL$ denote an invertible sheaf on $Y$ which has on every fiber
$Y_\fop$ the degree
$- \lceil \mu_{\rm min}(\shF_\eta (-m)) \rceil +1$.

Then
$ \mu_{\rm min}(\shF_\eta (-m) \otimes \shL_\eta)
= \mu_{\rm min}(\shF_\eta(-m))
- \lceil \mu_{\rm min}(\shF_\eta (-m))  \rceil +1 >0$
and hence $\shF_\eta (-m) \otimes \shL_\eta $ is ample on $Y_\eta$
due to \ref{amplekritnull}.
Therefore $\shF_\fop (-m) \otimes \shL_\fop $ is ample on
$Y_\fop$ for $p \gg 0$, since ampleness is an open property
(see \cite[Th\'{e}or\`{e}me 4.7.1]{EGAIII}
or \cite[Theorem 1.2.13]{lazarsfeld}; since we need here only
the generic open property, we can also use \ref{affineopen} together with
\cite[Th\'{e}or\`{e}me 4.5.2] {EGAII}).
This means again by \ref{amplekrit} that
$\bar{\mu}_{\rm min}( \shF_\fop(-m) \otimes \shL_\fop) >0$ or that
$$\bar{\mu}_{\rm min}( \shF_\fop(-m)) > - \deg \, (\shL_\fop)
= \lceil \mu_{\rm min}(\shF_\eta (-m))  \rceil - 1 \, .$$

This gives the first result.
The second statement follows by applying the first statement to $\shR(m)$,
\begin{eqnarray*}
\bar{\mu}_{\rm max} (\shF(-m)) &= &-\bar{\mu}_{\rm min} (\shR(m)) \cr
&<& -( \lceil \mu_{\rm min}(\shR_\eta (m))  \rceil -1) \cr
&=& \lfloor -\mu_{\rm min}(\shR_\eta (m)) \rfloor  +1 \cr
&=& \lfloor \mu_{\rm max}(\shF_\eta (-m)) \rfloor +1 \, .
\end{eqnarray*}
\qed

\section{Applications to tight closure: vanishing theorems}
\label{sectionvanishing}

The inclusion bound in Theorem \ref{maxin} and the exclusion bound in
Theorem \ref{minex} coincide if the sheaf of relations is semistable.
This gives a new class of vanishing type theorems in dimension two
and generalizes the vanishing theorem for para\-meter ideals,
see \cite{hunekesmithkodaira} and Corollary \ref{parametervanishing} below
(the name vanishing is due to the fact
that it is related to Kodaira Vanishing Theorem).
We give first the formulation in zero characteristic.

\begin{theorem}
\label{semistablevanishing}
Let the notation be as in {\rm \ref{forcingsit}}
and suppose that the characteristic of $K$ is zero.
Set $k= \lceil \frac{d_1 + \ldots +d_n}{n-1} \rceil $.
Suppose that the sheaf of relations $\shR(m)$ for ideal generators
$f_1, \ldots ,f_n$ is semistable. Then
$$(f_1, \ldots ,f_n)^* = (f_1, \ldots ,f_n) + R_{\geq k} \, .$$
\end{theorem}

\proof
Let $f_0 \in R$ be homogeneous of degree $m$.
Suppose first that $m \geq k$.
Then
$$m \geq \frac{d_1 + \ldots +d_n}{n-1} =\frac{\mu (\shF(0))}{\delta}
=\frac{\mu_{\max}(f_1, \ldots ,f_n)}{\delta} \, .$$
Hence the numerical condition in \ref{maxin} is fulfilled,
thus $f_0 \in (f_1, \ldots ,f_n)^*$.

Suppose now that $m<k$.
Then $m< \frac{d_1+ \ldots +d_n}{n-1}
=\frac{\mu_{\min}(f_1, \ldots, f_n)}{\delta}$
and \ref{minex} gives the result.
\qed

\medskip
Suppose that in the relative setting (\ref{relsituation})
the sheaf of relations is semistable in the generic point, so that
the vanishing theorem \ref{semistablevanishing} holds
in the generic point.
What can we say about the behavior in positive characteristic?
We know by \cite[\S 5]{miyaokachern} that
$\shF_\fop$ is semistable on an open non-empty subset
of $\Spec\,D$. However, for strongly semistable
we have to take into account the following problem of Miyaoka.

\begin{remark}
Miyaoka states in \cite[Problem 5.4]{miyaokachern} the following problem:
suppose that $C$ is a relative (smooth projective) curve over a
(say) $\ZZ$-algebra $D$ of finite type
and assume that a locally free sheaf $\shF$
is semistable in the generic fiber (characteristic zero).
Let $S$ be the set of points $P \in \Spec \, D$
of positive characteristic such that
$\shF | {C_P}$ is strongly semistable. Is $S$ dense in $\Spec\, D$?
\end{remark}

Therefore we may not expect that semistability in the generic point implies
a vanishing theorem for $p \gg 0$ without further conditions.
It implies however that the bounds are quite near to the expected number
$(d_1 + \ldots + d_n)/(n-1)$.

\begin{corollary}
\label{comparebound}
Suppose the situation of {\rm \ref{relsituation}} and suppose that
the sheaf of relations is semistable over the generic point.
Let $m =\deg(f_0)$. Then the following hold for $p \gg 0$.

\renewcommand{\labelenumi}{(\roman{enumi})}
\begin{enumerate}

\item
If $m \geq  \frac{d_1 +\ldots + d_n}{n-1} + \frac{1}{\delta} $, then
$f_0 \in (f_1, \ldots ,f_n)^*$.

\item
If $m \leq \frac{d_1 +\ldots + d_n}{n-1} - \frac{1}{\delta}$,
then
$f_0 \in (f_1, \ldots ,f_n)^*$ if and only if $f_0 \in (f_1, \ldots ,f_n)$.
\end{enumerate}
\end{corollary}
\proof
For (i) we have the estimates
\begin{eqnarray*}
\delta m & \geq &  \delta \frac{d_1+\ldots+d_n}{n-1} + 1 \cr
&\geq & \lfloor \delta \frac{d_1+\ldots+d_n}{n-1} \rfloor + 1 \cr
&= & \lfloor \mu_{\rm \max} (\shF_\eta) \rfloor +1 \cr
&> & \bar{\mu}_{\rm max}(\shF) \, ,
\end{eqnarray*}
where the last estimate follows from Proposition \ref{compare}.
The statement follows from Theorem \ref{maxin}.

(ii) follows by similar estimates from Theorem \ref{minex}.
\qed

\medskip
The previous corollary shows that the inclusion and exclusion bounds
are very near to $(d_1+ \ldots + d_n)/(n-1)$. If this number is not an
integer and if the degree of the curve is big enough, then we get also a
vanishing theorem from \ref{comparebound}.
In general however we get a vanishing theorem in positive characteristic only
for those points $\fop \in \Spec\, D$
for which the sheaf of relations $\shR_\fop$ is strongly semistable.

\begin{theorem}
\label{semistablevanishingp}
Suppose the situation and notation of \ref{relsituation}.
Set $k= \lceil \frac{d_1 + \ldots +d_n}{n-1} \rceil $.
Suppose that for the generic fiber
the sheaf of relations is semistable. Then for all $p \gg 0$
such that the corresponding sheaf of
relations is strongly semistable we have
$$(f_1, \ldots ,f_n)^* = (f_1, \ldots ,f_n) + R_{\geq k} \, .$$
\end{theorem}

\proof
The inclusion $\supseteq$ follows for every $\fop \in \Spec\, D$ such that
the sheaf $\shR_{\fop}(m)$
is strongly semistable on $Y_{\bar{\kappa}(\fop)}$
from \ref{maxin}.

For the inclusion $\subseteq $ we do not need the
condition strongly semistable.
From Theorem \ref{minex} we know that a single fixed element
$f_0$ with $f_0  \in (f_1, \ldots ,f_n)^*$
belongs also to the right hand side for $p \gg 0$,
but here we state the identity of the two ideals
for $p \gg 0$.
Let $I=(f_1, \ldots , f_n)$ in $S$. We may assume by shrinking
$D$ that $I=S \cap IS_{Q(D)}$.
Let $m <k$ and consider
$S_m /I_m \subseteq H^1(D(S_+),Rel(f_1, \ldots,f_n)_m)$.
We may assume that $S_m/I_m$ is a free $D$-module with a basis
induced by $h_j \in S_m$, $ 1 \leq j \leq t$.

For every field $Q(D) \subseteq L$
the sheaf $\shF_L(-m) = \shR_L(m)^\dual$ is ample on
$\Proj S_{L}$ due to \ref{amplekrit}, hence for every
$h = \sum \lambda_j h_j \neq 0 $
the extension $\shF'_L(-m)(h)$ is also ample
and the open subset $\PP(\shF'_{L,h})- \PP(\shF_L)$
is affine. This is then also true for the open subset
$U_L =D( (S_L)_+) \subseteq
\Spec \, S_L[T_1, \ldots ,T_n]/(f_1T_1 + \ldots +f_nT_n+h)$.

We introduce indeterminates $\Lambda_j, \, 1 \leq j \leq t,$
for the coefficients of an element $h= \sum_j \lambda_jh_j$
and consider the universal forcing algebra
$$C= S[T_i, \Lambda_j] /
(f_1T_1 + \ldots +f_nT_n + \Lambda_1h_1 + \ldots +\Lambda_th_t)$$
over $S[\Lambda_j]$ and over $D[\Lambda_j]$. Let $U=D(S_+) \subseteq \Spec\, C$.

We claim that (after shrinking $D$)
for every point $P \in \Spec \, D[\Lambda_1, \ldots , \Lambda_t]$,
$P \not\in V(\Lambda_1, \ldots , \Lambda_t)$,
the fiber $U_P$ is affine.
We show this by increasing inductively
the open subset where this statement holds.

For the quotient field $L=Q(D)(\Lambda_1, \ldots ,\Lambda_t)
=Q(D[\Lambda_1, \ldots , \Lambda_t])$
we know that $U_L$ is affine.
Therefore also $U \cap D(g)$ is affine (by remark \ref{affineopen}),
where $0 \neq g \in D[\Lambda_1, \ldots , \Lambda_t]$,
and we know then that the fiber $U_P$ is affine for every point
$P \in D(g) \in \Spec \, D[\Lambda_1, \ldots , \Lambda_t]$.
So we know that the claim is true
for a non-empty open subset.

For the induction step suppose that the claim is true for
the open subset $W \subseteq \Spec \, D[\Lambda_1 , \ldots, \Lambda_t]$.
By shrinking $D$ we may assume that all irreducible components
of the complement of $W$ dominate $\Spec \, D$.
Consider such an irreducible component $Z=V(\foq)$ of
maximal dimension
and suppose that $Z \neq V(\Lambda_1, \ldots , \Lambda_t)$.
The generic point $\zeta$ of $Z$ has characteristic zero
and the $\Lambda_j$ are not all zero in $\kappa(\zeta)$,
hence again $U_\zeta$ is affine and we find an open neighborhood
$\zeta \in D(g) \cap Z \subseteq Z$,
$g \in D[\Lambda_1, \ldots , \Lambda_t] $,
such that for every point $P \in D(g) \cap Z$
the fiber $U_P$ is affine.
Hence the claim is now true on a bigger open subset and the number of
components of the complement of maximal dimension has dropped.

So we see that the claim holds
eventually for $D(\Lambda_1, \ldots , \Lambda_t)$.
Now the claim
means in particular that for every prime ideal $\fop \in \Spec\, D$
and every linear combination
$h = \sum_j \lambda_j h_j \neq 0$, $ \lambda_j \in \kappa(\fop)$,
the corresponding
open subset $U_{\fop, h}$ is affine,
hence $h \not\in I^*_{\kappa(\fop)}$ for all $h \not\in I_{\kappa(\fop)}$.

This procedure can be done for every degree $0 \leq m < k$,
hence we find a sufficiently small $\Spec\, D$ such that
the statement holds for all $\fop \in \Spec\, D$.
\qed

\begin{remark}
The Theorems \ref{semistablevanishing}, \ref{comparebound}
and \ref{semistablevanishingp} indicate that the
number $\lceil \frac{d_1 + \ldots +d_n}{n-1} \rceil $
is the generic bound for the degree of an element
to belong to the tight closure,
see also Theorem \ref{genusbound}.
It is reasonable to guess that for $R$ of higher dimension the number
$\lceil \frac{ \dim \, R -1}{n-1} (d_1 + \ldots +d_n) \rceil $
should take over this part ($n \geq \dim \, R$).
\end{remark}

The following result of Huneke and Smith
\-- the vanishing theorem for parameters in dimension two \--
is an easy corollary of \ref{semistablevanishing} and
\ref{semistablevanishingp}, see \cite[Theorem 4.3]{hunekesmithkodaira}.

\begin{corollary}
\label{parametervanishing}
Let the notation be as in {\rm \ref{forcingsit}} and suppose that $n=2$,
so we are concerned with the tight closure of a parameter ideal.
Suppose that the characteristic of $K$ is $0$
or $p >  \frac{2}{\delta} (g-1)$.
Then
$$(f_1,f_2)^* = (f_1,f_2) + R_{\geq d_1+d_2} \, .$$
\end{corollary}
\proof
The statement for characteristic $0$
and characteristic $p \gg 0$
follows at once from \ref{semistablevanishing} and \ref{semistablevanishingp},
since the sheaf of relations $\shR$ is invertible, hence (strongly) semistable.

For the precise statement in positive characteristic we need
the ampleness criterion of Hartshorne-Mumford for bundles of rank two,
see \ref{amplecritranktwo}.
The inclusion $ \supseteq $ is true in any characteristic
by \ref{maxin}, since $\shR(m)= \O_Y(m-d_1-d_2)$.
So suppose that $f_0 \in (f_1,f_2)^*$, but
$f_0 \not\in (f_1,f_2)$ and $m= \deg \, (f_0) < d_1+d_2$.
The corresponding forcing class $c \in H^1(Y, \shR(m))$
gives the forcing sequence
$$ 0   \lra \O_Y \lra \shF'(-m) \lra \shF(-m)  \lra 0 \, .$$
Now $\deg \, (\shF'(-m)) = (d_1+d_2 -m) \delta \geq \delta$,
hence $\deg \, (\shF'(-m)) > \frac{2}{p} (g-1) $.
This gives the first condition in the criterion \ref{amplecritranktwo}.
For the second condition, let
$\shF'(-m) \ra \shM \ra 0$ be an invertible quotient sheaf.
If $\deg \, (\shM) < 0$ or if $ \deg \, (\shM) = 0$ and $\shM \neq \O_Y$,
then the composed mapping $\O_Y \ra \shM$ is zero and $\shM$ is a quotient
of the ample invertible sheaf $\shF(-m)$, which is not possible.
If $\shM= \O_Y$, then again the composed mapping
$\O_Y \ra \O_Y$ is not zero,
hence it is an isomorphism and the sequence splits,
which contradicts the assumption $c \neq 0$.
\qed

\begin{remark}
If $Y$ is a smooth plane curve given by a polynomial $F$ of degree $\delta$,
then $g=(\delta -1) (\delta -2) /2$ and the condition in
\ref{parametervanishing} is that
$p > \frac{2}{\delta} (\frac{( \delta -1) (\delta -2)}{2} -1) = \delta -3 $.

The condition in \cite[Theorem 4.3]{hunekesmithkodaira}
for \ref{parametervanishing} to hold in positive characteristic is that
$p$ exceeds the degree of each of a set of
homogeneous generators for the module of $k$-derivations of $R$.
\end{remark}

\section{Examples for $n=3$}
\label{examples}

In this last section we discuss the case where $n=3$. This means that
we have three homogeneous elements $f_1,f_2,f_3 \in R$ of degree $d_1,d_2,d_3$,
where $R$ is our two-dimensional, normal standard-graded
domain over an algebraically closed field.
Then the sheaf of relations has rank two.
If it is decomposable, then we are in the situation of
Theorem \ref{splitting} and the tight closure
$(f_1,f_2,f_3)^*$ is easy to compute, see the example \ref{indepexample}
below. So we shall concentrate on the indecomposable case.
We will restrict to characteristic $0$;
for positive charcteristic we may obtain with the help of
\ref{compare} results which are slightly worse.

\begin{theorem}
\label{genusbound}
Let $K$ denote an algebraically closed field of characteristic $0$
and let $R$ denote a normal
standard-graded $K$-domain of dimension two. Let $Y= \Proj\, R$
denote the corresponding smooth projective curve of genus $g$
and let $\delta$ denote the degree of $\O_Y(1)$.
Let $f_1,f_2,f_3$ denote $R_+$-primary homogeneous
elements of degree $d_1,d_2,d_3$
and suppose that the sheaf of relations $\shR(m)$
is indecomposable on $Y$. Then
$$\mu_{\rm max} (f_1,f_2,f_3) \leq \delta \frac{d_1+d_2+d_3}{2} + g-1 \, $$
and
$$ \mu_{\rm min} (f_1,f_2,f_3) \geq \delta \frac{d_1+d_2+d_3}{2} - g+1 \, .$$
\end{theorem}

\proof
The sheaf $\shF(-m)=\shR^\dual(m)$ has rank two,
hence $\PP(\shF)$ is a ruled surface over
$Y$. Since $\shF(-m)$ is supposed to be indecomposable,
we know by \cite[Theorem V.2.12]{haralg} that
the $e$-invariant of $\PP(\shF)$ is $e \leq 2g-2$.
If $\shF(-m)$ is semistable, then the result is true,
since then $\mu_{\rm max}(\shF(0)) = \mu(\shF(0)) = \delta (d_1+d_2+d_3)/2$
and $g \geq 1$, since $\shF(0)$ is indecomposable.

So suppose that $\shF(0)$ is not semistable, set $\shF= \shF(0)$ and let
$0 \ra \shL \ra \shF \ra \shQ \ra 0$ be a short exact sequence
such that $\deg \, (\shL) = \mu(\shL) =\mu_{\rm max} (\shF)$.
Then $\shF \otimes \shL^{-1}$ is normalized,
since it has a section $\neq 0$ and
$H^0( Y, \shF \otimes \shL^{-1} \otimes \shM) = 0$
for every invertible sheaf $\shM$ of  $\deg \, (\shM) <0$.
For else there would exist a non-trivial morphism
$ \shL \otimes \shM^{-1} \ra \shF$ which contradicts
$\mu_{\rm max}(\shF) = \mu(\shL)$. Hence
$$-e= \deg \, (\shF \otimes  \shL^{-1}) =2 \mu(\shF \otimes  \shL^{-1})=
2 \mu (\shF) +2 \mu(\shL^{-1}) = 2 \mu (\shF) -2 \mu(\shL)$$
and therefore
$$\mu_{\rm max}(\shF)= \mu(\shL) = \mu(\shF) + \frac{e}{2}
\leq \delta \frac{d_1+d_2+d_3}{2} + \frac{2g-2}{2} \, .$$
The bound for $\mu_{\rm min}(f_1,f_2,f_3)$
follows from $\mu_{\rm min}( \shF ) = 2 \mu(\shF) - \mu_{\rm max} (\shF)$.
\qed

\begin{corollary}
\label{genusboundtight}
Let $K$ denote an algebraically closed field of characteristic $0$
and let $R$ denote a normal standard-graded
$K$-domain of dimension two. Let $Y= \Proj\, R$
denote the corresponding smooth projective curve of genus $g$
and let $\delta$ denote the degree of $\O_Y(1)$.
Let $f_1,f_2,f_3$ denote $R_+$-primary homogeneous
elements of degree $d_1,d_2,d_3$
and suppose that the sheaf of relations $\shR$ is
indecomposable on $Y$. Then the following statements hold.

\renewcommand{\labelenumi}{(\roman{enumi})}
\begin{enumerate}

\item
For $m \geq \frac{d_1+d_2+d_3}{2} + \frac{g-1}{\delta} $
we have the inclusion $R_m \subseteq (f_1,f_2,f_3)^*$.

\item
For $m < \frac{d_1+d_2+d_3}{2} - \frac{g-1}{\delta}$
we have $(f_1,f_2,f_3)^* \cap R_m = (f_1,f_2,f_3) \cap R_m$.
\end{enumerate}
\end{corollary}

\proof
The first statement follows from Theorem \ref{genusbound}
and from Theorem \ref{maxin}.
The second statement follows from \ref{genusbound} and \ref{minex}
\qed

\begin{corollary}
\label{degreebound}
Let $K$ denote an algebraically closed field
of charcteristic $0$ and let $F \in K[x,y,z]$ denote a
homogeneous polynomial of degree $\delta$
such that $R =K[x,y,z]/(F)$ is a normal domain.
Let $f_1,f_2,f_3 \in R$ denote $R_+$-primary homogeneous
elements of degree $d_1,d_2,d_3$.
Suppose that the sheaf of relations $\shR$ is indecomposable
on the curve $Y= \Proj \, R$. Then the following hold.

\renewcommand{\labelenumi}{(\roman{enumi})}
\begin{enumerate}

\item
$R_m \subseteq (f_1,f_2,f_3)^*$
for $m \geq \frac{d_1+d_2+d_3}{2} + \frac{\delta -3}{2} $.

\item
For $m < \frac{d_1+d_2+d_3}{2} - \frac{\delta +3}{2}$
we have $(f_1,f_2,f_3)^* \cap R_m = (f_1,f_2,f_3) \cap R_m$.
\end{enumerate}

\end{corollary}

\proof
This follows from Corollary \ref{genusboundtight} taking into account that
$g=(d-1)(d-2)/2$.
\qed

\medskip
We also have the following result in positive characteristic.

\begin{corollary}
\label{genusboundpositive}
Suppose the relative setting {\rm \ref{relsituation}} and suppose that $n=3$.
Suppose that  $\shR(m)$ is indecomposable
on $Y_{\bar{\eta}}$, where $\eta $ is the generic point.
Then for $p \gg 0$ we have
$$ \bar{\mu}_{\rm max}(f_1,f_2,f_3)
< \lfloor \delta \frac{d_1+d_2+d_3}{2} \rfloor + g $$
and 
$$ \bar{\mu}_{\rm min}(f_1,f_2,f_3)
> \lfloor \delta \frac{d_1+d_2+d_3}{2} \rfloor - g $$
\end{corollary}
\proof
We have from \ref{compare} and \ref{genusbound} the estimates
\begin{eqnarray*}
\bar{\mu}_{\rm max} (f_1,f_2,f_3) &< &
\lfloor \mu_{\rm max}(\shF_\eta (0)) \rfloor  +1 \cr
& \leq & \lfloor  \delta \frac{d_1+d_2+d_3}{2} + g-1 \rfloor +1 \cr
&=& \lfloor  \delta \frac{d_1+d_2+d_3}{2} \rfloor +g 
\end{eqnarray*}
\qed

\begin{example}
\label{150example}
Let $F \in K[x,y,z]$ denote a homogeneous form of degree $\delta =5$
such that $R=K[x,y,z]/(F)$ is a normal domain.
Suppose that $f_1,f_2,f_3 \in R$ are primary
homogeneous elements of degree $100$ and such that its sheaf of relations
is indecomposable. Then the bounds of Smith in
\cite{smithgraded} give $R_{\geq 200} \subseteq (f_1,f_2,f_3)^*$ on the one
hand and, on the other hand, that an element of degree $ m \leq 100$ belongs
to $(f_1,f_2,f_3)^*$ if and only if it belongs to the ideal itself.

Our bounds in \ref{degreebound}
give $R_{\geq 151} \subseteq (f_1,f_2,f_3)^*$ and
that an element of degree $\leq  148$ belongs
to $(f_1,f_2,f_3)^*$ if and only if it belongs to the ideal itself
(for zero characteristic).
So only the degrees $149$ and $150$ are not covered by \ref{degreebound}.
If furthermore the sheaf of relations is semistable,
then \ref{semistablevanishing} gives a complete numerical answer.
\end{example}

\begin{remark}
It seems not so easy to establish the semistability property
for the relation bundle on a given curve of genus $\geq 2$.
There exist however restriction theorems saying that this property holds
on a generic curve of sufficiently high degree under the condition
that the bundle is (defined and) semistable on the projective plane.
We will come back to this kind of questions in another paper.
\end{remark}

Our last example gives a negative answer to a question
of Craig Huneke asked at the MSRI (September 2002).
The example gives an ideal which is generated by
$*$-independent elements (meaning that none of them is contained
in the tight closure of the others, see \cite{vraciu} for this notion),
but it does not hold a vanishing theorem for it, that is
there does not exist a common inclusion and exclusion
bound for tight closure.

\begin{example}
\label{indepexample}
Consider the ideal $I=(x^4,xy,y^2)$ in the Fermat cubic $x^3+y^3+z^3=0$.
These ideal generators come from
the regular polynomial ring $K[x,y] \subset K[x,y,z]/(x^3+y^3+z^3)=R$.
From this it follows that they are $*$-independent in $R$
and that the relation bundle must split, and
in fact
$$\shR(x^4,xy,y^2)(5) = \O_Y(0) \oplus \O_Y(2) \,$$
where
$Y= \Proj \, R$ is the corresponding elliptic curve.
Let $h$ denote a homogeneous element of degree $m$, given rise
to a forcing class
$$c \in H^1(Y, \shR(m))= H^1(Y,\O_Y(m-5)) \oplus H^1(Y,\O_Y(m-3)) \, . $$
From the numerical criterion in Theorem \ref{splitting}
it is easy to deduce the following.

For $m \geq 5$ we have $R_m \subset I^*$.

For $m=4$ an element $h$ belongs to $I^*$ only if it belongs to $I$.

For $m=3$ all possibilities occur.
We find
$ yz^2 \in (x^4,xy,y^2)^*$, but not in the ideal itself,
and $ xz^2 \not\in I^*$.
\end{example}


\end{document}